\documentclass[12pt]{article}
\usepackage{epic,eepic,epsf,epsfig}
\usepackage{amsmath}
\usepackage{amssymb}
\usepackage{amsbsy}

\usepackage{amsfonts}
\newtheorem{theorem}{Theorem}[section]%
\newtheorem{cor}[theorem]{Corollary}%
\newtheorem{remark}[theorem]{Remark}%

\def\f{\noindent}

\newcommand{\qed}{\mbox{\raisebox{0.7ex}{\fbox{}}} \vspace{4truemm}}
\def\demo{\f {\bf Proof.}\hskip10pt}

\begin{document}

\baselineskip 16pt

\title{ \vspace{-1.2cm}
Finite groups in which every maximal invariant subgroup of order divisible by $p$ is nilpotent
\thanks{\scriptsize This research was supported in part by Shandong Provincial Natural Science Foundation, China (ZR2017MA022)
and NSFC (11761079).
\newline
 \hspace*{0.5cm} \scriptsize $^{\ast\ast}$Corresponding
  author.
\newline
       \hspace*{0.5cm} \scriptsize{E-mail addresses:}
       shijt2005@163.com\,(J. Shi),\,shanmj2023@s.ytu.edu.cn\,(M. Shan),\,xufj2023@s.ytu.edu.cn\,(F. Xu).}}

\author{Jiangtao Shi\,$^{\ast\ast}$,\,Mengjiao Shan,\,Fanjie Xu\\
\\
{\small School of Mathematics and Information Sciences, Yantai University, Yantai 264005, China}}

\date{ }

\maketitle \vspace{-.8cm}

\begin{abstract}
Let $A$ and $G$ be finite groups such that $A$ acts coprimely on $G$ by automorphisms. For any fixed prime divisor $p$ of $|G|$, we provide a
complete characterization of the structure of a group $G$ in which every maximal $A$-invariant subgroup of order divisible by $p$ is nilpotent.

\medskip \f {\bf Keywords:} maximal invariant subgroup; Sylow subgroup; nilpotent group
\end{abstract}

\section{Introduction}

All groups are considered to be finite. Recall that a group $G$ in which every maximal subgroup is nilpotent is called a Schmidt group, and
R$\rm\acute{e}$dei {\rm\cite{redei1956}} provided a complete characterization of Schmidt groups. As a generalization, Deng, Meng and Lu
{\rm\cite[Theorem 3.1]{meng}} gave a complete description of the structure of a group of even order in which every maximal subgroup of
even order is nilpotent. Furthermore, Shi and Tian {\rm\cite[Theorem 1.1]{shitian}} obtained a complete classification of a group
$G$ in which every maximal subgroup of order divisible by $p$ is nilpotent for any fixed prime divisor $p$ of $|G|$.

As another generalization of Schmidt groups, Beltr$\rm\acute{a}$n and Shao {\rm\cite[Theorem A]{beltran}} had the following result under
the hypothesis of coprime action of groups.

\begin{theorem} {\rm\cite[Theorem A]{beltran}}\ \ Let $G$ and $A$ be groups of coprime orders and assume that $A$ acts on $G$
by automorphisms. If every maximal $A$-invariant subgroup of $G$ is nilpotent but $G$ is not,
then $G$ is solvable and $|G|=p^aq^b$ for two distinct primes $p$ and $q$, and $G$ has a normal
$A$-invariant Sylow subgroup.
\end{theorem}

Furthermore, considering any fixed prime divisor $p$ of $G|$, Beltr$\rm\acute{a}$n and Shao {\rm\cite[Theorem B]{beltran23}} obtained the
following result.

\begin{theorem} {\rm\cite[Theorem B]{beltran23}}\ \ Suppose that a group $A$ acts coprimely on a group $G$ and let $p$ be
a prime divisor of the order of $G$. If every maximal $A$-invariant subgroup of
$G$ whose order is divisible by $p$ is nilpotent, then $G$ is soluble.
\end{theorem}

In this paper, our main goal is to provide a complete characterization of a group $G$ in which every maximal $A$-invariant
subgroup of order divisible by $p$ is nilpotent. We have the following result, the proof of which is given is Section~\ref{s2}.

\begin{theorem}\ \ \label{th1} Let $A$ and $G$ be groups such that $A$ acts coprimely on $G$ by
automorphisms. Then every maximal $A$-invariant
subgroup $G$ of order divisible by $p$ is nilpotent for any fixed prime divisor $p$ of $|G|$ if and only if one of the following statements holds.

$(1)$ $G$ is nilpotent;

$(2)$ $G=Q\rtimes P$, where $Q$ is a normal Sylow $q$-subgroup of $G$ for $q\neq p$ and $Q$ has an $A$-invariant proper subgroup $Q_0$
such that $Q_0$ is normal in $G$ and $Q_0\times P$ is a unique nilpotent maximal $A$-invariant subgroup of $G$ containing $P$, $P$ is an
$A$-invariant Sylow $p$-subgroup of $G$ having a unique maximal $A$-invariant subgroup $P_0$ and $[Q,P_0]=1$.

$(3)$ $G=P\rtimes Q$, where $P$ is a normal Sylow $p$-subgroup of $G$ and $P$ has an $A$-invariant proper subgroup $P_0$ such
that $P_0$ is normal in $G$ and $P_0\times Q$ is a unique nilpotent maximal $A$-invariant subgroup of $G$ containing $Q$, $Q$ is an $A$-invariant
Sylow $q$-subgroup of $G$ having a unique maximal $A$-invariant subgroup $Q_0$ and $[P,Q_0]=1$, $q\neq p$.

$(4)$ $G=P\times(Q\rtimes R)$, where $P$ is a normal Sylow $p$-subgroup of $G$, $Q$ is a normal Sylow $q$-subgroup of $G$ and $Q$ has an $A$-invariant proper subgroup $Q_0$ of $Q$ such that $Q_0$ is normal in $Q\rtimes R$ and $Q_0\times R$ is a unique nilpotent maximal
$A$-invariant subgroup of $Q\rtimes R$ containing $R$, $R$ is an $A$-invariant Sylow $r$-subgroup of $G$ having a unique maximal $A$-invariant subgroup
$R_0$ and $[Q,R_0]=1$, $Q\rtimes R$ is a unique non-nilpotent maximal $A$-invariant subgroup of $G$, $p$, $q$ and $r$ are distinct primes.
\end{theorem}

\begin{remark} {\rm It is obvious that a group $G$ having a unique maximal subgroup is a cyclic group of prime-power order. However, a group $G$ having a unique maximal $A$-invariant subgroup might not be a cyclic group. For example, let $G=Q_8$ and $A$ be a subgroup of ${\rm Aut}(G)$ of order 3. It is is easy to see that $G$ has a unique maximal $A$-invariant subgroup of order 2 but $G$ is non-cyclic.}
\end{remark}

\begin{remark} {\rm Note that a group $G$ is a minimal non-nilpotent group if and only if $G=P\rtimes Q$, where $P$ is a normal Sylow subgroup of $G$ having a proper subgroup $P_0$ such that $P_0$ is normal in $G$ and $P_0\times Q$ is a unique nilpotent maximal subgroup of $G$ containing $Q$, $Q$ is a cyclic Sylow subgroup of $G$. Hence {\rm\cite[Theorem 1.1]{shitian}} is a direct corollary of Theorem~\ref{th1} when $A=1$.}
\end{remark}

The following corollary is a direct consequence of Theorem~\ref{th1}.

\begin{cor}\ \ Let $A$ and $G$ be groups such that $A$ acts coprimely on $G$ by automorphisms. Then $G$ is non-$p$-nilpotent and every maximal $A$-invariant
subgroup $G$ of order divisible by $p$ is nilpotent for any fixed prime divisor $p$ of $|G|$ if and only if $G=P\rtimes Q$, where $P$ is a normal Sylow $p$-subgroup of $G$ and $P$ has an $A$-invariant proper subgroup $P_0$ such
that $P_0$ is normal in $G$ and $P_0\times Q$ is a unique nilpotent maximal $A$-invariant subgroup of $G$ containing $Q$, $Q$ is an $A$-invariant
Sylow $q$-subgroup of $G$ having a unique maximal $A$-invariant subgroup $Q_0$ and $[P,Q_0]=1$, $q\neq p$.
\end{cor}

\section{Proof of Theorem~\ref{th1}}\label{s2}

\demo The sufficiency part is evident, we only need to prove the sufficiency part.

We first prove that $G$ has normal Sylow subgroups.

Otherwise, assume that $G$ has no normal Sylow subgroups. We divide our analyzes into two cases.

By the hypothesis, every maximal $A$-invariant subgroup of $G$ of is $p$-nilpotent.

Case (1): Suppose that $G$ is not $p$-nilpotent. Then the Sylow $p$-subgroup of $G$ is normal by {\rm\cite[Corollary 4]{meng2020}}, a contradiction.

Case (2): Suppose that $G$ is $p$-nilpotent. Let $P$ be an $A$-invariant Sylow $p$-subgroup of $G$ and $M$ be a normal $p$-complement of $P$ in $G$. Then $G=M\rtimes P$. Note that $M$ is an $A$-invariant normal Hall-subgroup of $G$.
By our assumption, $P$ is not normal in $G$.

For $M$, when $M$ is nilpotent, it is obvious that $G$ has normal Sylow subgroups, a contradiction. Therefore, $M$ is non-nilpotent. By the hypothesis, it is easy to see that $M$ must be a non-nilpotent maximal $A$-invariant subgroup of $G$. Moreover, $M$ must be a unique non-nilpotent maximal $A$-invariant subgroup of $G$. Let $S$ be an $A$-invariant Sylow subgroup of $M$. By our assumption, $S$ is not normal in $G$. Then $N_G(S)<G$. It is obvious that $N_G(S)\nleq M$. Therefore, there exists a nilpotent maximal $A$-invariant subgroup $L$ of $G$ such that $N_G(S)\leq L$.
Since $G=ML$, one has $|P|\mid|L|$. Let $P_{01}$ be an $A$-invariant Sylow $p$-subgroup of $L$, then $N_G(P_{01})=L$ and $|S|\mid|L|$. Since $M$ is non-nilpotent,
there exists an $A$-invariant Sylow subgroup $T$ of $M$ such that $(|S|,|T|)=1$. Arguing as above, there exits a nilpotent maximal $A$-invariant subgroup $N$ of $G$
such that $N_G(T)\leq N$ and $|P|\mid|N|$. Let $P_{02}$ be an $A$-invariant Sylow $p$-subgroup of $N$, one has $N_G(P_{02})=N$ and $|T|\mid|N|$. Note that $P_{01}$
and $P_{02}$ are conjugate in $G$, which implies that $N_G(P_{01})$ and $N_G(P_{02})$ are conjugate in $G$. It follows that $|L|=|N|$ and then the order of every $A$-invariant Sylow subgroup of $G$ divides $|L|$, a contradiction.

Hence our assumption is not true and then $G$ has normal Sylow subgroups.

Let $P_1,\,P_2,\,\cdots,\,P_s$ be all normal Sylow subgroups of $G$. Assume $E=P_1\times P_2\times\cdots\times P_s$. Since $G$ is solvable by {\rm\cite[Theorem B]{beltran23}}, there exists an $A$-invariant subgroup $F$ of $G$ such that $G=E\rtimes F$.

When $p\mid|F|$, one has that $F$ is nilpotent by the hypothesis. Assume $F=Q_1\times Q_2\times\cdots\times Q_t$, where $Q_1=P,\,Q_2,\,\cdots,\,Q_t$ are
$A$-invariant Sylow subgroups of $F$, $t\geq 1$.

If $t>1$, then $E\rtimes Q_1=E\rtimes P$ is nilpotent, which implies that $P$ is normal in $G$,
a contradiction. Therefore, $t=1$. It follows that $G=E\rtimes P$.

If $s>1$, then $(P_1\times\cdots\times P_{i-1}\times P_{i+1}\times\cdots\times P_s)\rtimes P$
is nilpotent for every $1\leq i\leq s$, which implies that $[P_i,P]=1$ for every $1\leq i\leq s$ and then $G$ is nilpotent, a contradiction.
Therefore, $s=1$.

Assume $Q=P_1$, then $G=Q\rtimes P$. It is easy to see that every maximal $A$-invariant subgroup of $G$ either contains $Q$ or contains an $A$-invariant
Sylow $p$-subgroup of $G$. Then every maximal $A$-invariant subgroup of $G$ is nilpotent by the hypothesis.

For $P$. If $P$ has two distinct $A$-invariant maximal subgroups $P_{11}$ and $P_{12}$. Since both $P_{11}$ and $P_{12}$ are normal in $P$, one has $P=P_{11}P_{12}$. By the hypothesis, both $Q\rtimes P_{11}$
and $Q\rtimes P_{12}$ are nilpotent, which implies that $[Q,P_{11}]=1$ and $[Q,P_{12}]=1$. It follows that $[Q,P]=1$ and then $G$ is nilpotent, a contradiction. Therefore, $P$ has a unique maximal $A$-invariant subgroup $P_0$ and $[Q,P_0]=1$.

Let $B$ be a maximal $A$-invariant subgroup of $G$ such that $P\leq B$. Then $B=(B\cap Q)\rtimes P$. One has that $B$ is nilpotent by the hypothesis. It is obvious that $B\cap Q<Q$, which implies that $N_G(B\cap Q)>B$.
It follows that $B\cap Q$ is normal in $G$ by the maximality of $B$. Let $Q_0=B\cap Q$. If there exists a maximal $A$-invariant subgroup $B_0$ of $G$ containing $P$ and $B_0\neq Q_0\times P$, then $B_0\cap Q\neq Q_0$.
Note that $Q_0\neq 1$ and $B_0\cap Q\neq 1$. Moreover, both $Q_0$ and $B_0\cap Q$ are normal in $G$. It is obvious that $Q_0\nleq B_0$. Then $G=Q_0B_0=Q_0((B_0\cap Q)P)$. It follows that $Q=Q_0(B_0\cap Q)$. Since
$[Q_0,P]=1$ and $[B_0\cap Q,P]=1$, one has $[Q,P]=1$, which implies that $G$ is nilpotent, a contradiction. Therefore, $Q_0\times P$ is a unique nilpotent maximal $A$-invariant subgroup of $G$ containing $P$.

When $p\nmid|F|$, then $p\mid|E|$. For $E=P_1\times P_2\times\cdots\times P_s$, assume $P=P_1$. Arguing as above, one has $s=1$ or $s=2$.

If $s=1$, then $G=P\rtimes F$.

When $F$ is nilpotent, it is easy to see that $F$ must be an $A$-invariant Sylow subgroup of $G$. Assume $F=Q$. Arguing as above, one has that $Q$ has a unique maximal $A$-invariant subgroup $Q_0$ and $[P,Q_0]=1$, and $P$ has an $A$-invariant proper subgroup $P_0$ such that $P_0$ is normal in $G$ and $P_0\times Q$ is a unique nilpotent maximal $A$-invariant subgroup of $G$ containing $Q$.

When $F$ is non-nilpotent, then every maximal $A$-invariant subgroup of $F$ is nilpotent by the hypothesis. By {\rm\cite[Theorem A]{beltran}}, assume $F=Q\rtimes R$, where $Q$ and $R$ are $A$-invariant Sylow subgroups of $F$. Then $G=P\rtimes(Q\rtimes R)$. By the hypothesis,
both $P\rtimes Q$ and $P\rtimes R$ are nilpotent. It follows that $G=P\times(Q\rtimes R)$. Note that every maximal $A$-invariant subgroup of $G$
either contains $P$ or contains some conjugate of $Q\rtimes R$. Since $G$ is non-nilpotent, $Q\rtimes R$ is non-nilpotent. It follows that $Q\rtimes R$ is a
non-nilpotent maximal $A$-invariant subgroup of $G$. Moreover, $Q\rtimes R$ is a unique non-nilpotent maximal $A$-invariant subgroup of $G$. Arguing as
above, one has that $R$ has a unique maximal $A$-invariant subgroup $R_0$ and $[Q,R_0]=1$, $Q$ has an $A$-invariant proper subgroup $Q_0$ such that $Q_0$ is normal in $Q\rtimes R$ and $Q_0\times R$ is a unique nilpotent maximal $A$-invariant subgroup of $Q\rtimes R$ containing $R$.

If $s=2$, assume $P=P_1$ and $Q=P_2$. One has $G=(P\times Q)\rtimes F$. Obviously, $F$ must be an $A$-invariant Sylow subgroup of $G$ by the hypothesis. Then one can easily get that $G$ belongs to the above case when $s=1$ and $F$ is non-nilpotent.\hfill\qed


\begin{thebibliography}{99}
\bibitem{shitian} J. Shi and Y. Tian, On finite groups in which every maximal subgroup of order divisible by $p$ is nilpotent (or abelian),  Rendiconti del Seminario Matematico della Universit$\rm\grave{a}$ di Padova,
DOI: 10.4171/RSMUP/159.
\end{thebibliography}
\end{document}